\documentclass[12pt,a4paper]{article}
\usepackage{localeng,hyperref}
\begin{document}
\init
\title{Comparing angles in Euclid's \emph{Elements}}
\author{Alexander Shen\thanks{LIRMM, Univ Montpellier, CNRS, Montpellier, France, ORCID~0000-0001-8605-7734\\ \texttt{sasha.shen@gmail.com}, \texttt{alexander.shen@lirmm.fr}}}
\date{}
\maketitle

\begin{abstract}
The exposition in Euclid's Elements contains an obvious gap (seemingly unnoticed by most commentators):  he often compares not just angles, but \emph{groups} of angles, and at the same time he avoids summing angles (and considering angles greater than $\pi$), and does not say what such a comparison of groups could mean. We discuss the problem and suggest a possible interpretation that could make Euclid's exposition consistent.
\end{abstract}

This note records my attempts to understand the basic part of Euclid's Elements --- not from the historian's viewpoint\footnote{I am not a historian and cannot read the original Greek text; therefore I used only English and Russian translations --- mainly~\cite{english,english-heath,russian}.} but from the naive reader's viewpoint. Imagine a student (with some modern mathematics background) who is preparing to pass Euclid's exam and wants to know what should be said during the exam and what should not.

There are many places in Euclid's text that look unclear. But some of them are harmless. For example, when we see the axiom ``A point is that of which is no part'', we can just note that this axiom is never used, so it is enough to learn this sentence by heart (just in case) and relax. There are some other locally unclear places.  For example,  how should we explain that two different lines cannot have two different common points? (There is some special axiom for that in some editions, but it is unclear what Euclid would say about that.) But let us focus on some recurring theme: the comparison of quantities.

The first pages of Elements contain ``common notions'' (called sometimes ``axioms'') that say, e.g.,  ``And if equal things are added to equal things then the wholes are equal''~\cite[Book~1, page~7]{english}. However, in Book~5 Euclid explains the theory of proportions using a different language (``magnitudes'' in~\cite{english}, and indeed it seems that the Greek word is also different), but it is not clear whether \emph{things} from Book 1 and \emph{magnitudes} from Book 5 are the same objects. Anyway, in both cases compared objects include line segments (called ``lines''), angles and figures (figures are ``equal'' for Euclid if they have the same area; they do not need to be congruent).

Let us concentrate on angles. It seems that Euclid considers only angles that are (in modern terminology) strictly smaller than $\pi$. Angles of size $0$ and $\pi$ do not exist (``And a plane angle is the inclination of the lines to one another when two lines in a plane meet one another, and \emph{are not lying in a straight line}'', Definition~8 from Book~1), and it seems that Euclid never considers angles greater than~$\pi$, either. The next definition shows that Euclid also considers angles whose sides are not lines (line segments) but curves: ``and when the lines containing the angle are straight then the angle is called rectilinear''. However, almost all the time Euclid uses only rectilinear angles. (Still in Book 3, proposition 16~\cite[p.~87]{english} Euclid claims that the angle between a semicircle and the diameter connecting its endpoints is greater that any acute angle.)

What about comparing angles? Only one \emph{common notion} speaks about ``greater'' relation (``And the whole [is] greater than the part''); all other speak only about equal/unequal things. Still some properties of the \emph{greater} relation are used freely in the proofs. For example, Proposition~16 (Book~1) assumes that if one angle is equal to another one that is greater than the third one, then the first one is greater than the third one~\cite[p.~21, lines 7--8]{english}.

The next natural question: do we assume that angles are linearly ordered, i.e., every two angles $\alpha$ and $\beta$ are either equal, or one of them is greater?  The proof of Proposition~25~\cite[p.~28]{english}  seems to use this property: ``angle $BAC$ is not equal to $EDF$. Neither, indeed, is $BAC$ less than $EDF$. $\langle\ldots\rangle$ Thus, $BAC$ is greater than $EDF$.'' But there is no axiom about linear ordering. May be one should prove it by applying one angle to another  (one side to one side), and then saying that obviously one of the angles is inside the other and therefore is smaller (because \emph{the whole is greater than the part})? The superposition of figures seems to be a legal argument for Euclid, it is used in Proposition~4 of Book~1 when proving that two triangles with the same pair of sides and angle between them are equal. This proof uses  \emph{common notion}~4: \emph{things coinciding with one another are equal to one another}. However, in the next Proposition 5 saying that for isosceles triangles, the angles at the base are equal to one another, Euclid does not use the obvious argument (rotate the triangle and apply it to itself exchanging the sides) and considers instead some additional construction. (Why?)

Note also that if we allow angles with arbitrary curves as sides, we should probably allow them to be incomparable because it is hard to see what kind of linear order can be established on them. (Some later commentators, e.g. Proclus, seem to agree with that.) But let us restrict ourselves to rectilinear angles, since we have a bigger problem ahead: is it allowed to add angles?  The 5th postulate says ``if a straight-line falling across two (other) straight-lines makes internal angles on the same side (of itself whose sum is) less than two right-angles, then the two (other) straight-lines, being produced to infinity, meet on that side''~\cite[p.~7]{english}. The words in parentheses are added by translator (without explanation)\footnote{Another English translation~\cite[p.~202]{english-heath} says simply ``internal angles on the same side less than two right angles'', without the word ``sum''.}. Russian translation~\cite{russian} says just ``angles smaller that two right angles''\footnote{\rus{углы, меньшие двух прямых}}, and the translator's note explains: ``We would say now `angles whose \emph{sum} is smaller than two right angles', but we keep here and below the Euclid's formulation''\footnote{\rus{Мы сказали бы <<углы, \emph{в сумме} меньшие двух прямых>>. Здесь и далее мы сохраняем евклидов способ изложения.}} The same problem appears in \emph{common notion}~2: ``if equal things are added to equal things then the wholes are equal'', says the English translation. Russian translation comments ``Euclid says `\gre{τὰ ὅλα}', and the correct translation is ``whole'', not ``sums''; Euclid does not considers adding of quantities and the resulting sums''\footnote{\rus{<<У Евклида `\gre{τὰ ὅλα}', что вполне точно переводится `целые', а не `суммы'; Евклид не мыслит сложения величин и получаемых после сложения сумм>>}.}

But if Euclid does not consider sums of angles, what does he mean by ``two angles smaller than two right angles''? What kind of objects does he compare? And if he does consider sums of angles, what kind of objects these sums are? What is the sum of two right angles? It is not an angle, since the angle of size $\pi$ does not exist for Euclid. Moreover, in Book 4, Proposition 3~\cite[p.~112]{english} Euclid speaks about sums of four right angles: ``the (sum of the) four angles of quadrilateral $AMBK$ is equal to four right angles''. Maybe, the sum of two right angle is a halfplane, and the sum of four right angles is the entire plane, so these sums are geometric figures (though not angles), and we can compare them even if they are not angles\footnote{Would Euclid agree that the sum of angles of a hexagon is equal to eight right angles?}? 

Probably, Euclid would not agree with this interpretation either. Indeed, Proposition 13 from Book 1 says that two adjacent angles $ABC$ and $DBA$ are together equal to two right angles (whatever it means). 
\begin{center}
\includegraphics[width=0.4\textwidth]{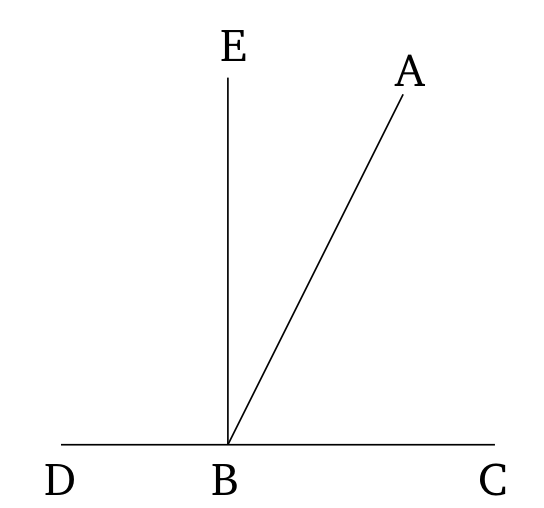}
\end{center}
Euclid draws $EB$ that forms two right angles $EBD$ and $EBC$. Now it seems that the sum of $DBA$ and $ABC$ as a geometric figure is the same as the sum of two right angles $EBD$ and $EBC$ (both sums are the same half-plane), so the proposition is proved. But Euclid does not think in this way. Instead, he says that $CBE$ is equal to $CBA$ and $ABE$, and continues: ``let $EBD$ have been added to both. Thus, the (sum of the angles) $CBE$ and $EBD$ is equal to the (sum of the three) angles $CBA$, $ABE$ and $EBD$''. Then he shows in a similar way that the (sum of the angles) $DBA$ and $ABC$ is also equal to the (sum of the three angles) $CBA$, $ABE$, and $EBD$, and it remains to use the transitivity of equality.

How can all these arguments be interpreted in a consistent way? What are the ``things'' compared according to common notions, and what properties of these ``things'' are used?

Unfortunately, it seems that commented Euclid's editions (including all mentioned in the bibliography) do not take this question seriously. Often they provide a lot of historical context (who said what, where and when), or try to comment on further developments (that are now obsolete), but the basic question about \emph{how we could interpret Euclid's text in a consistent way} remains mostly unanswered.\footnote{For example, there is an adaptation~\cite{english-casey} that its author descibes as follows: ``a work which, while giving the unrivalled original in all its integrity, would also contain the modern conceptions and developments of the portion of Geometry over which the Elements extend''~\cite[Preface]{english-casey}. The author tries to preserve the sequence and numbering of propositions but changes the definitions (and their numbering). He says explicitly ``The angle by joining two or more angles together is called their sum'' in the definitions (which Euclid never does) but later (following Euclid) uses the word ``sum'' for the sum of four right angles, without any explanations. Also he tacitly omits the parts where Euclid speaks about angles whose sides are curves.} It seems to me that the only (more or less) consistent interpretation of what Euclid says about angles goes as follows.

\begin{itemize}

\item ``Angles'' are geometric angles (pairs of rays) between $0$ and $\pi$.

\item ``Things'' that are considered in \emph{common notions} are not angles but \emph{finite multisets of angles} (in other words, commutative and associative formal sums of angles).

\item There is a linear order on those ``things''.

\item Every singleton is greater than the empty multiset.

\item The order is preserved when some angle is added (as a new element) to both sides of a comparison.

\item If a ray splits an angle $\alpha$ into two angles $\beta$ and $\gamma$, then $\{\alpha\}$ and $\{\beta,\gamma\}$ are equal multisets. 

\item Angles that coincide when applied to each other are equal.

\end{itemize}

If we understand ``things'' like that\footnote{This approach somehow resembles local groups. In~\cite{english-joyce} some remarks can be interpreted in a similar way: ``That sum being mentioned is a straight angle, which is not to be considered as an angle according to Euclid. It is a formal sum equal to two right angles. In other propositions formal sums of four right angles occur. These and larger formal sums are not angles themselves, merely sums of angles. Only if an angle sum is less than two right angles can it be identified with a single angle'' (comments to Proposition 13 of Book 1). Still in the comments to Common Notions the author writes ``These common notions, sometimes called axioms, refer to magnitudes of one kind. The various kinds of magnitudes that occur in the Elements include lines, angles, plane figures, and solid figures'', not mentioning ``formal sums'', so again there is no consistent picture.}, the Euclid's reasoning becomes understandable (and one can understand why he proves the theorem about adjacent angles in this seemingly strange way). Still the results from Book 5 about proportions needs some adjustments (if the``magnitudes'' are understood in the same way, as multisets). Say, Proposition 8~\cite[p.~139]{english} explicitly uses the fact that the given magnitude is a singleton, not an arbitrary multiset. But this gaps can be filled in the proposed framework (the problems with Archimedes' axiom remain, but this is another story).

This setting looks artificial, but one could illustrate it by an analogy: the balance can be used not only to compare two weights, but to compare two \emph{sets} of weights. Then we assume that we can add the same weight to both sides not changing the balance, and that we can split a weight in two without influencing the balance. Maybe, Euclid or some commentators already used this analogy to explain the theory of angles?

\emph{Acknowledgements}. The author is grateful to Andrei Rodin, Alexander Shtern, Andrei Shchetnikov and the participants of the Kolmogorov complexity seminar for their comments and explanations.

\end{document}